\documentclass[final]{amsart} 

\usepackage{amsmath, amssymb, latexsym, amscd}
\usepackage{exscale, cite, eps fig, graphics}

\usepackage[colorlinks]{hyperref}
\usepackage[notref, notcite]{showkeys}

\renewcommand{\geq}{\geqslant}
\renewcommand{\leq}{\leqslant}

\renewcommand{\l}{\langle}
\renewcommand{\r}{\rangle}

\renewcommand{\c}{c_{\rm ww}}
\renewcommand{\L}{\mathcal{L}}
\newcommand{\M}{\mathcal{M}}

\renewcommand{\k}{\kappa}

\newcommand{\R}{\mathbb{R}}

\newcommand{\C}{\mathbb{C}}

\newtheorem{theorem}{Theorem}[section]
\newtheorem{lemma}[theorem]{Lemma}

\newtheorem*{main-theorem}{Main Theorem}
\newtheorem*{remark*}{Remark}

\numberwithin{equation}{section}

\title[Modulational instability in the FDCH equation]
{Modulational instability in \\ the full-dispersion Camassa-Holm equation}

\author[Hur]{Vera~Mikyoung~Hur}
\address{Department of Mathematics, University of Illinois at Urbana-Champaign, Urbana, IL 61801 USA}
\email{verahur@math.uiuc.edu}

\author[Pandey]{Ashish~K.~Pandey}
\email{akpande2@illinois.edu}  

\date{\today}

\begin{document}

\maketitle

\begin{abstract}
We determine the stability and instability of a sufficiently small and periodic traveling wave to long wavelength perturbations, for a nonlinear dispersive equation which extends a Camassa-Holm equation to include all the dispersion of water waves and the Whitham equation to include nonlinearities of medium amplitude waves. In the absence of the effects of surface tension, the result qualitatively agrees with the Benjamin-Feir instability of a Stokes wave. In the presence of the effects of surface tension, it qualitatively agrees with those from formal asymptotic expansions of the physical problem and it improves upon that for the Whitham equation, correctly predicting the limit of strong surface tension. We discuss the modulational stability and instability in the Camassa-Holm equation and related models.
\end{abstract}

\section{Introduction}\label{sec:intro}

In the 1960s, Whitham (see \cite{Whitham}, for instance) proposed
\begin{equation}\label{E:whitham0}
\eta_t+c_{\rm ww}(\sqrt{\beta}|\partial_x|)\eta_x+(3\sqrt{(1+\alpha\eta)}-3)\eta_x=0,
\end{equation}
to argue for wave breaking in shallow water. That is, the solution remains bounded but its slope becomes unbounded in finite time. Here $t\in\mathbb{R}$ is proportional to elapsed time, and $x\in\mathbb{R}$ is the spatial variable in the primary direction of wave propagation; $\eta=\eta(x,t)$ is the fluid surface displacement from the undisturbed depth, 
\[
\alpha=\frac{\text{a typical amplitude}}{\text{the undisturbed fluid depth}}\quad\text{and}\quad
\beta=\frac{(\text{the undisturbed fluid depth})^2}{(\text{a typical wavelength})^2}.
\]
Moreover, $c_{\rm ww}(|\partial_x|)$ is a Fourier multiplier operator, defined as 
\begin{equation}\label{def:c1}
\widehat{c_{\rm ww}(|\partial_x|)f}(\k)=\sqrt{\frac{\tanh\k}{\k}}\widehat{f}(\k).
\end{equation}
Note that $\c(\k)$ means the phase speed in the linear theory of water waves. 
For small amplitude waves satisfying $\alpha\ll1$, we may expand the nonlinearity of \eqref{E:whitham0} up to terms of order $\alpha$ to arrive at
\begin{equation}\label{E:whitham1}
\eta_t+c_{\rm ww}(\sqrt{\beta}|\partial_x|)\eta_x+\frac32\alpha\eta\eta_x=0.
\end{equation}
For relatively shallow water or, equivalently, relatively long waves satisfying $\beta\ll1$, we may expand the right side of \eqref{def:c1} up to terms of order $\beta$ to find
\[
c_{\rm ww}(\sqrt{\beta}k)=1-\frac16\beta k^2+O(\beta^2).
\]
Therefore, for small amplitude and long waves satisfying $\alpha=O(\beta)$ and $\beta\ll1$, we arrive at the famous Korteweg-de Vries equation
\begin{equation}\label{E:KdV}
\eta_t+\eta_x+\frac16\beta\eta_{xxx}+\frac32\alpha\eta\eta_x=0.
\end{equation}
As a matter of fact, for well-prepared initial data, the solutions of the Whitham equation and the Korteweg-de Vries equation differ from those of the water wave problem merely by higher order terms over the relevant time scale; see \cite{Lannes}, for instance, for details. But \eqref{E:whitham1} and \eqref{def:c1} offer improvements over \eqref{E:KdV} for short waves. Whitham conjectured wave breaking for \eqref{E:whitham1} and \eqref{def:c1}. One of the authors \cite{Hur-breaking} recently proved this. In stark contrast, no solutions of \eqref{E:KdV} break. 

Moreover, Johnson and one of the authors \cite{HJ2} showed that a sufficiently small and $2\pi/\k$ periodic traveling wave of the Whitham equation be spectrally unstable to long wavelength perturbations, provided that $\k>1.145\dots$. In other words, \eqref{E:whitham1} (or \eqref{E:whitham0}) and \eqref{def:c1} predict the Benjamin-Feir instability of a Stokes wave; see  \cite{BF, Benjamin1967, Whitham1967} and \cite{BM1995}, for instance. In contrast, periodic traveling waves of the Korteweg-de Vries equation are all modulationally stable. By the way, under the assumption that $\eta_t+\eta_x$ is small, we may modify \eqref{E:KdV} to arrive at the Benjamin-Bona-Mahony equation
\begin{equation}\label{E:BBM}
\eta_t+\eta_x-\frac16\beta\eta_{xxt}+\frac32\alpha\eta\eta_x=0.
\end{equation}
It agrees with \eqref{E:KdV} for long waves but is preferable for short waves. Note that the phase speed for \eqref{E:BBM} is bounded for all frequencies. The authors \cite{HP1} showed that a sufficiently small and $2\pi/\k$ periodic traveling wave of \eqref{E:BBM} be modulationally unstable if $\k>\sqrt{3}$. Hence the Benjamin-Bona-Mahony equation seems to predict the Benjamin-Feir instability of a Stokes wave. But the instability mechanism is different from that in the Whitham equation or the water wave problem; see \cite{HP1} for details. 

Furthermore, in the presence of the effects of surface tension, Johnson and one of the authors \cite{HJ3} determined the modulational stability and instability of a sufficiently small and periodic traveling wave of \eqref{E:whitham1} and 
\begin{equation}\label{def:cT1}
\widehat{c_{\rm ww}(|\partial_x|;T)f}(\k)=\sqrt{(1+T\k^2)\frac{\tanh\k}{\k}}\widehat{f}(\k),
\end{equation}
where $T$ is the coefficient of surface tension. The result agrees by and large with those in \cite{Kawahara, DR}, for instance, from formal asymptotic expansions of the physical problem. But it fails to predict the limit of ``strong surface tension." Perhaps, this is not surprising because \eqref{E:whitham1} neglects higher order nonlinearities of the water wave problem. It is interesting to find an equation, which predicts the modulational stability and instability of a gravity capillary wave. This is the subject of investigation here. 

By the way, the authors \cite{HP2} recently extended the Whitham equation to include bidirectional propagation, and they showed that the ``full-dispersion shallow water equations" correctly predict capillary effects on the Benjamin-Feir instability. But the modulation calculation is very lengthy and tedious. Here we seek higher order nonlinearities suitable for unidirectional propagation. Such a model is likely to be Hamiltonian and is potentially useful for other purposes.

As a matter of fact, for medium amplitude and long waves satisfying $\alpha=O(\sqrt{\beta})$ and $\beta\ll1$, the Camassa-Holm equations for the fluid surface displacement
\begin{equation}\label{E:CH1}
\eta_t+\eta_x+\beta(a\eta_{xxx}+b\eta_{xxt})+\frac32\alpha\eta\eta_x-\frac38\alpha^2\eta^2\eta_x+\frac{3}{16}\alpha^3\eta^3\eta_x=-\alpha\beta(c\eta\eta_{xxx}+d\eta_x\eta_{xx})\hspace*{-20pt}
\end{equation}
and for the average horizontal velocity
\begin{equation}\label{E:CHu1}
u_t+u_x+\beta(au_{xxx}+bu_{xxt})+\frac32\alpha uu_x=-\alpha\beta(cuu_{xxx}+du_xu_{xx}),
\end{equation}
where 
\[
0\leq a\leq \frac16, \qquad b=a-\frac16,\qquad c=\frac32a+\frac16,\quad\text{and}\quad d=\frac92a+\frac{5}{24},
\]
extend the Korteweg-de Vries equation to include higher order nonlinearities, and they approximate the physical problem; see \cite{Lannes}, for instance, for details. In the case of $a=1/12$, \eqref{E:CH1} reads
\[
\eta_t+\eta_x+\frac{1}{12}\beta(\eta_{xxx}-\eta_{xxt})+\frac32\alpha\eta\eta_x-\frac38\alpha^2\eta^2\eta_x+\frac{3}{16}\alpha^3\eta^3\eta_x=-\frac{7}{24}\alpha\beta(\eta\eta_{xxx}+2\eta_x\eta_{xx}),
\]
which is particularly interesting because it predicts wave breaking; see \cite{Lannes} and references therein. Note that 
\[
\frac{3\alpha\eta}{1+\sqrt{1+\alpha\eta}}=\frac32\alpha\eta-\frac38\alpha^2\eta^2+\frac{3}{16}\alpha^3\eta^3+O(\alpha^4).
\]
Lannes~\cite{Lannes} combined the dispersion relation of water waves and a Camassa-Holm equation, to propose the {\em full-dispersion Camassa-Holm} (FDCH) equation for the fluid surface displacement
\begin{equation}\label{E:FDCH1}
\eta_t+c_{\rm ww}(\sqrt{\beta}|\partial_x|)\eta_x+\frac{3\alpha\eta}{1+\sqrt{1+\alpha\eta}}\eta_x
=-\alpha\beta\Big(\frac{5}{12}\eta\eta_{xxx}+\frac{23}{24}\eta_x\eta_{xx}\Big),
\end{equation}
where $\c(|\partial_x|)$ is in \eqref{def:c1}, or \eqref{def:cT1} in the presence of the effects of surface tension. For relatively long waves satisfying $\beta\ll1$, \eqref{E:FDCH1} and \eqref{def:c1} agree with \eqref{E:CH1}, where $a=1/6$, up to terms of order $\beta$. But, including all the dispersion of water waves, \eqref{E:FDCH1} and \eqref{def:c1} may offer an improvement over \eqref{E:CH1} for short waves. For small amplitude waves satisfying $\alpha\ll1$, \eqref{E:FDCH1} agrees with \eqref{E:whitham1} up to terms of order $\alpha$. But, including higher order nonlinearities, \eqref{E:FDCH1} may offer an improvement over \eqref{E:whitham1} for medium amplitude waves. 
For the average horizontal velocity, we may combine \eqref{def:c1}, or \eqref{def:cT1} in the presence of the effects of surface tension, and \eqref{E:CHu1} to introduce
\begin{equation}\label{E:FDCHu1}
u_t+c_{\rm ww}(\sqrt{\beta}|\partial_x|)u_x+\frac32\alpha uu_x
=-\alpha\beta\Big(\frac{5}{12}uu_{xxx}+\frac{23}{24}u_xu_{xx}\Big).
\end{equation}

We follow along the same line as the arguments in \cite{HJ2, HJ3, HP1} (see also \cite{BHJ}) and investigate the modulational stability and instability in the FDCH equation. A main difference lies in that the nonlinearities of \eqref{E:FDCH1} involve higher order derivatives and, hence, a periodic traveling wave is not a priori smooth. We examine the mapping properties of various operators to construct a smooth solution. 

In the absence of the effects of surface tension, we show that a sufficiently small and $2\pi/\k$ periodic traveling wave of \eqref{E:FDCH1} and \eqref{def:c1} is spectrally unstable to long wavelength perturbations, provided that 
\[
\k>1.420\dots,
\] 
and stable to square integrable perturbations otherwise. The result qualitatively agrees with the Benjamin-Feir instability of a Stokes wave (see \cite{BF, Benjamin1967, Whitham1967}, for instance) and that for the Whitham equation (see \cite{HJ2}). The critical wave number compares reasonably well with $1.363\dots$ in the Benjmain-Feir instability. Including the effects of surface tension, in the $\k$ and $\k\sqrt{T}$ plane, we determine the regions of modulational stability and instability for a sufficiently small and periodic traveling wave of \eqref{E:FDCH1} and \eqref{def:cT1}; see Figure~\ref{f:ST} for details. The result qualitatively agrees with those in \cite{Kawahara, DR}, for instance, from formal asymptotic expansions of the physical problem, and it improves upon that in \cite{HJ3} for the Whitham equation. In particular, the limit of $\k(T)\sqrt{T}\to1.283\dots$ as $T\to \infty$, where $\k(T)$ is a critical wave number, whereas the limit is unbounded for the Whitham equation (see \cite{HJ3}). 

Moreover, we show that a sufficiently small and $2\pi/\k$ periodic traveling wave of \eqref{E:CH1} is modulationally unstable if $\k>6$. To the best of the authors' knowledge, this is new. The Camassa-Holm equation seems to predict the Benjamin-Feir instability of a Stokes wave. But the instability mechanism is different from that in \eqref{E:FDCH1} and \eqref{def:c1}, or the water wave problem. One may use the Evans function and other ODE methods to determine the modulational stability and instability for all amplitudes. This is an interesting direction of future research. The result herein indicates that the stability and instability depend on the carrier wave. 

In the absence of the effects of surface tension, we show that a sufficiently small and $2\pi/\k$ periodic traveling wave of \eqref{E:FDCHu1} and \eqref{def:c1} is modulationally unstable if $\k$ is greater than a critical value, similarly to the Benjamin-Feir instability. But, in the presence of the effects of surface tension, the modulational stability and instability in \eqref{E:FDCHu1} and \eqref{def:cT1} qualitatively agree with that in the Whitham equation (see \cite{HJ3}). In particular, it fails to predict the limit of strong surface tension. Therefore, we learn that the higher power nonlinearities of \eqref{E:FDCH1} improve the result, not the higher derivative nonlinearities. 

It is interesting to explain breaking, peaking, and other phenomena of water waves in \eqref{E:FDCH1} and \eqref{def:c1} (or \eqref{def:cT1}).

\subsubsection*{Notation}

The notation to be used is mostly standard, but worth briefly reviewing. Let $\mathbb{T}$ denote the unit circle in $\mathbb{C}$. We identify functions over $\mathbb{T}$ with $2\pi$                                                                                                                                                                                                    periodic functions over $\mathbb{R}$ via $f(e^{iz})=F(z)$ and, for simplicity of notation, we write $f(z)$ rather than $f(e^{iz})$. For $p$ in the range $[1,\infty]$, let $L^p(\mathbb{T})$ consist of real or complex valued, Lebesgue measurable, and $2\pi$ periodic functions over $\mathbb{R}$ such that 
\[
\|f\|_{L^p(\mathbb{T})}:=\Big(\frac{1}{2\pi}\int^\pi_{-\pi}|f(z)|^p~dz\Big)^{1/p}<\infty
\qquad\text{if\quad$p<\infty$},
\]
and $\|f\|_{L^\infty(\mathbb{T})}:=\text{ess\,sup}_{-\pi<z\leq\pi}|f(z)|<\infty$ if $p=\infty$. Let $H^1(\mathbb{T})$ consist of $L^2(\mathbb{T})$ functions whose derivatives are in $L^2(\mathbb{T})$. Let $H^\infty(\mathbb{T})=\bigcap_{k=0}^\infty H^k(\mathbb{T})$. 

For $f \in L^1(\mathbb{T})$, the Fourier series of $f$ is defined by
\[
\sum_{n\in \mathbb{Z}}\widehat{f}(n)e^{inz}, \qquad\text{where}\quad
\widehat{f}(n)=\frac{1}{2\pi}\int^\pi_{-\pi}f(z)e^{-inz}~dz.
\]
If $f\in L^2(\mathbb{T})$ then its Fourier series converges to $f$ pointwise almost everywhere. We define the $L^2(\mathbb{T})$ inner product as 
\begin{equation}\label{def:prod1}
\langle f_1,f_2\rangle_{L^2(\mathbb{T})}=\frac{1}{2\pi}\int^\pi_{-\pi} f_1(z)f_2^*(z)~dz
=\sum_{n\in\mathbb{Z}}\widehat{f_1}(n)\widehat{f_2}^*(n).
\end{equation}

\section{Sufficiently small and periodic traveling waves}\label{sec:periodic}

We determine periodic traveling waves of the FDCH equation, after normalization of parameters, 
\begin{equation}\label{E:FDCH}
\eta_t+c_{\rm ww}(|\partial_x|;T)\eta_x+\frac{3\eta}{1+\sqrt{1+\eta}}\eta_x
=-\Big(\frac{5}{12}\eta\eta_{xxx}+\frac{23}{24}\eta_x\eta_{xx}\Big),
\end{equation}
where $\c(|\partial_x|;T)$ is in \eqref{def:cT1}, and we calculate their small amplitude expansion.

\subsection*{Properties of $\c(\cdot\,;T)$}

For any $T\geq 0$, $\c(\cdot\,;T)$ is even and real analytic, and $\c(0;T)=1$. Note that $\c(|\partial_x|;0):H^s(\mathbb{R})\to H^{s+1/2}(\mathbb{R})$ for any $s\in\mathbb{R}$, and for $T>0$, $\c(|\partial_x|;T):H^{s+1/2}(\mathbb{R})\to H^{s}(\mathbb{R})$; see \cite{HJ2,HJ3}, for instance, for details.

\begin{figure}[h]
(a)~~\includegraphics[scale=0.3]{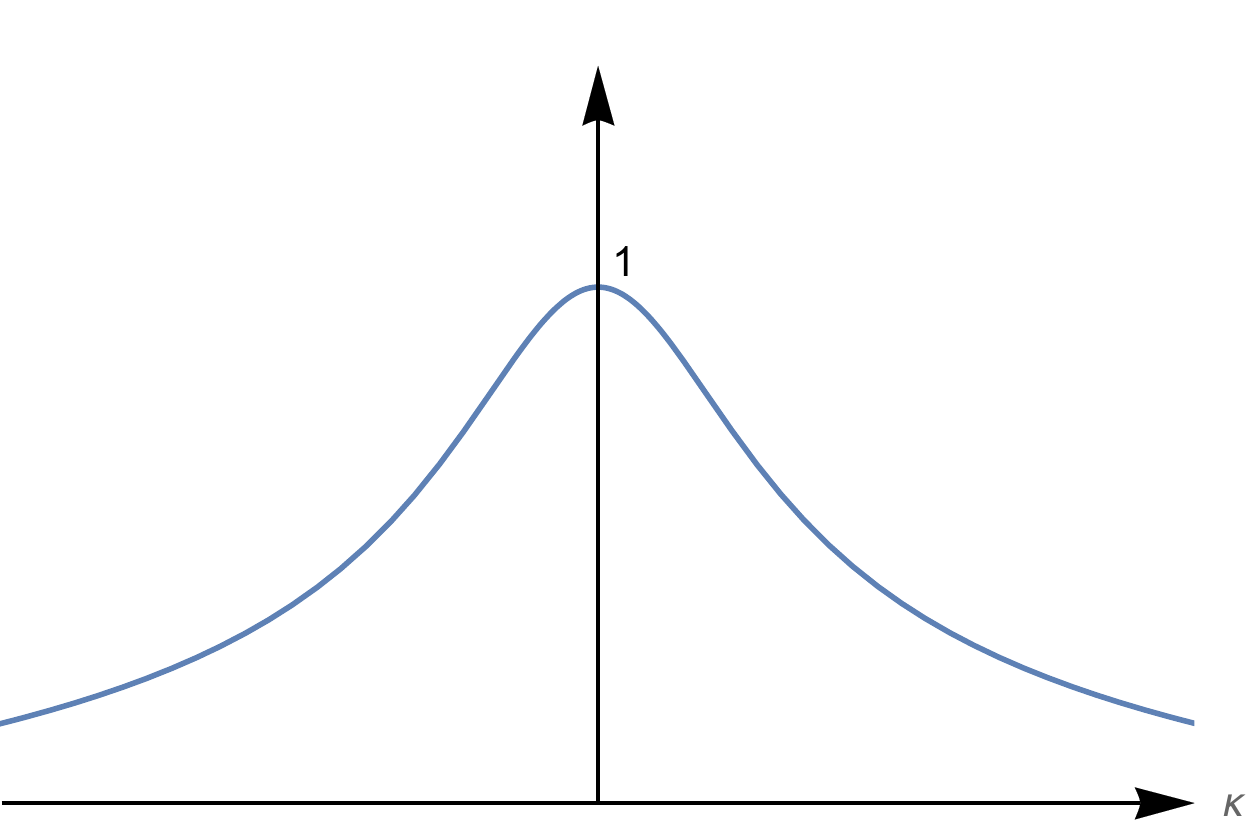}\quad(b)~~\includegraphics[scale=0.3]{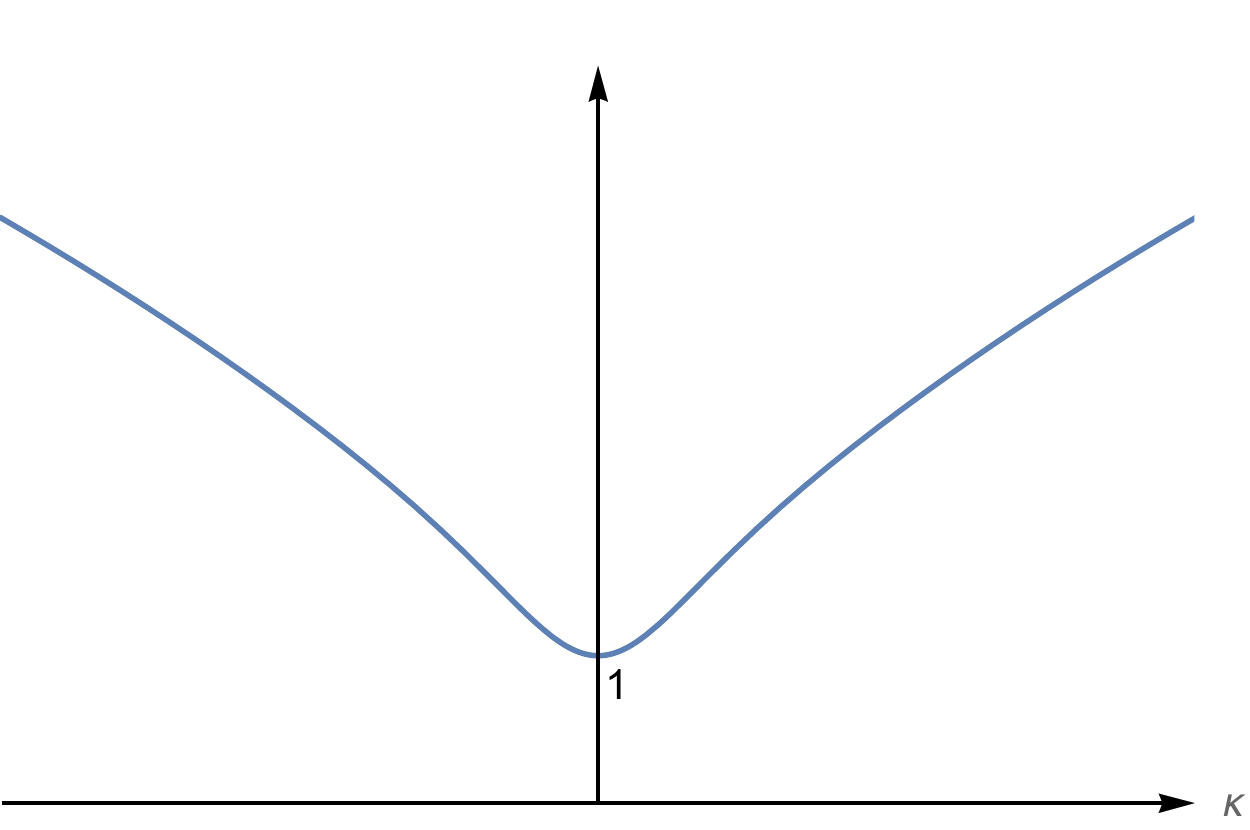}\quad(c)~~\includegraphics[scale=0.3]{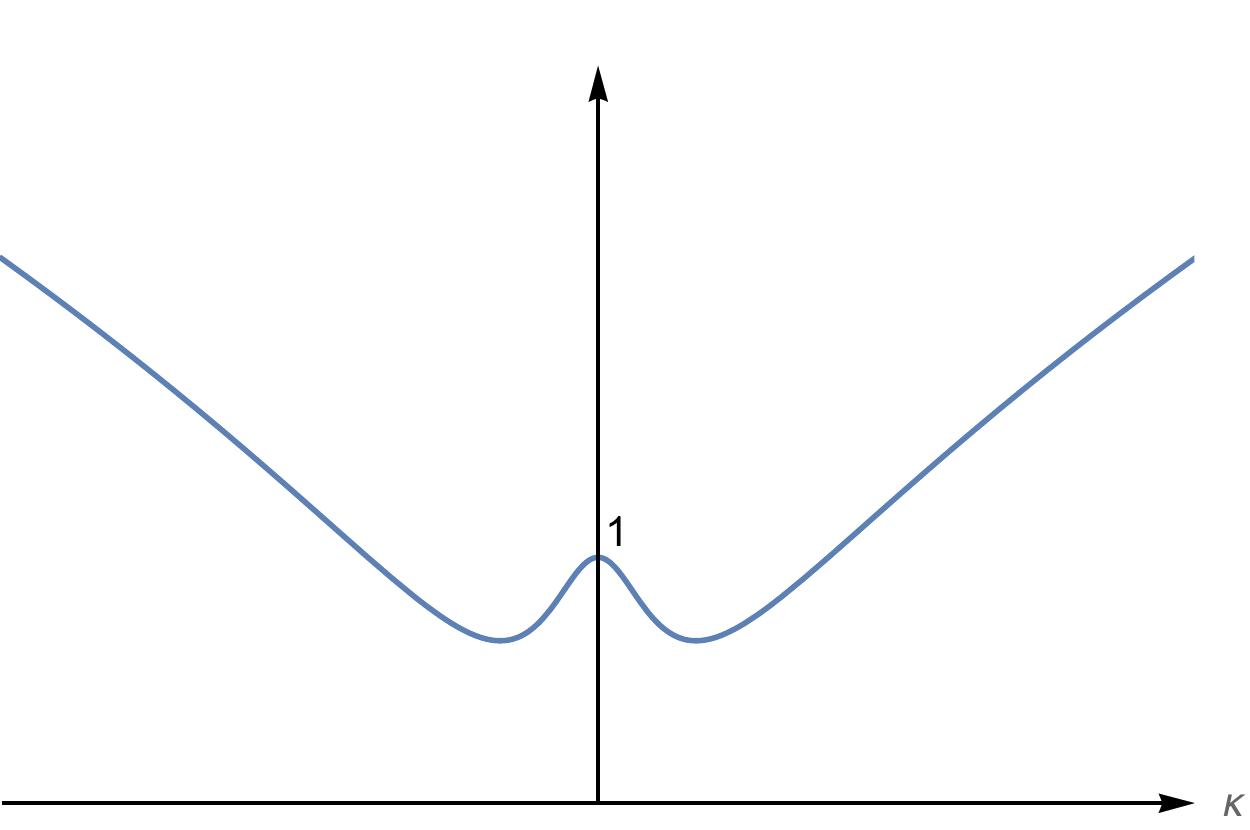}
\caption{Schematic plots of $\c(\cdot\,;T)$ for (a)~$T=0$, (b)~$T\geq1/3$, and (c)~$0<T<1/3$.}\label{fig:cT}
\end{figure}

Note that $\c(\cdot;0)$ decreases to zero monotonically away from the origin. For $T\geq1/3$, $\c(\cdot\,;T)$ increases monotonically and unboundedly away from the origin. For $0<T<1/3$, on the other hand, $\c'(0;T)=0$, $\c''(0;T)<0$ and $\c(\k;T)\to\infty$ as $\k\to\infty$. Hence $\c(\cdot\,;T)$ possesses a unique minimum over the interval $(0,\infty)$; see Figure~\ref{fig:cT}. 

By a traveling wave of \eqref{E:FDCH} and \eqref{def:cT1}, we mean a solution of the form $\eta(x,t)=\eta(x-ct)$ 
for some $c>0$, the wave speed, where $\eta$ satisfies by quadrature
\[
(\c(|\partial_x|;T)-c-3)\eta+2(1+\eta)^{3/2}-2+\frac{5}{12}\eta \eta_{xx}+\frac{13}{48} \eta_x^2 = (1-c)^2 b
\]
for some $b\in\mathbb{R}$. We seek a periodic traveling wave of \eqref{E:FDCH} and \eqref{def:cT1}. That is, $\eta$ is a $2\pi$ periodic function of $z:=\k x$ for some $\k>0$, the wave number, and it satisfies
\begin{equation}\label{E:pFDCH}
(\c(\k|\partial_z|;T) -c -3)\eta +2(1+\eta)^{3/2}-2+\frac{5}{12}\k^2\eta\eta_{zz}+\frac{13}{48}\k^2\eta_z^2 =(1-c)^2b.
\end{equation}
Note that
\begin{equation}\label{E:c-Hk}
\c(\k|\partial_z|;0): H^s(\mathbb{T})\to H^{s+1/2}(\mathbb{T}),\quad\text{and}\quad
\c(\k|\partial_z|;T): H^{s+1/2}(\mathbb{T})\to H^{s}(\mathbb{T})
\end{equation}
for $T>0$, for any $\k>0$ and $s\in\mathbb{R}$.
Note that
\begin{equation}\label{def:ck}
\c(\k|\partial_z|;T)e^{inz}=\c(n\k;T)e^{inz}\quad \text{for $n\in\mathbb{Z}$}.
\end{equation} 

Note that \eqref{E:pFDCH} remains invariant under 
\begin{equation}\label{E:invariance}
z\mapsto z+z_0\quad\text{and}\quad z\mapsto -z
\end{equation}
for any $z_0\in \mathbb{R}$. Hence we may assume that $\eta$ is even. But \eqref{E:pFDCH} does not possess scaling invariance. Hence we may not a priori assume that $\k=1$. Rather, the (in)stability result herein depends on the carrier wave number. Moreover, \eqref{E:pFDCH} does not possess Galilean invariance. Hence we may not a priori assume that $b=0$. Rather, we exploit the variation of \eqref{E:pFDCH} in the $b$ variable in the instability proof. To compare, the Whitham equation for periodic traveling waves possesses Galilean invariance; see \cite{HJ3}, for instance.

We follow along the same line as the arguments in \cite{HJ2,HJ3,HP1}, for instance, to construct periodic traveling waves of \eqref{E:FDCH} and \eqref{def:cT1}. A main difference lies in the lack of a priori smoothness of solutions of \eqref{E:pFDCH}. We examine the mapping properties of various operators to construct smooth solutions.

For any $T\geq0$ and an integer $k\geq0$, let 
\[
F:H^{k+2}(\mathbb{T})\times\mathbb{R}_+\times \mathbb{R}\times\mathbb{R}_+\to H^k(\mathbb{T})
\] 
denote
\begin{multline}\label{def:F}
F(\eta,c;b,\k,T)=(\c(\k|\partial_z|;T) -c -3)\eta +2(1+\eta)^{3/2}-2\\
\quad+\frac{5}{12} \k^2 \eta \eta_{zz}+\frac{13}{48} \k^2 \eta_z^2 -(1-c)^2 b.
\end{multline}
It is well defined by \eqref{E:c-Hk} and a Sobolev inequality. We seek a solution $\eta\in H^{k+2}(\mathbb{T})$, $c>0$, and $b\in \mathbb{R}$ of 
\begin{equation}\label{E:F=0}
F(\eta,c;b, \k,T)=0.
\end{equation}
Since $k$ is arbitrary, $\eta\in H^\infty(\mathbb{T})$. Note that $F$ is invariant under \eqref{E:invariance}. Hence we may assume that $\eta$ is even. 

For any $T\geq0$, and $c>0$, $b\in\mathbb{R}$, $\k>0$, note that
\begin{multline*}
F_\eta(\eta,c; b, \k,T)\zeta=\Big(\c(\k|\partial_z|;T)-c-3+3(1+\eta)^{1/2}\\
+\k^2\Big(\frac{5}{12}(\eta_{zz}+\eta \partial_z^2)+\frac{13}{24}\eta_z \partial_z\Big)\Big)\zeta:
H^{k+2}(\mathbb{T})\to H^k(\mathbb{T})
\end{multline*}
is continuous by \eqref{E:c-Hk} and a Sobolev inequality. Here a subscript means Fr\'echet differentiation. Moreover, for any $T\geq0$, and $\eta\in H^{k+2}(\mathbb{T})$, $\k>0$, $b\in\mathbb{R}$, note that $F_c(\eta;\k,c,b)=-\eta+2(1-c)b:\mathbb{R}\to H^k(\mathbb{T})$ is continuous. Since $F_b(\eta;\k,c,b)=-(1-c)^2$ and 
\[
F_\k(\eta;\k,c,b):=\c'(\k |\partial_z|;T) \eta + \frac{5}{6} \k \eta\eta_{zz}+\frac{13}{24} \k \eta_z^2
\]
are continuous likewise, $F$ depends continuously differentiably on its arguments. Furthermore, since the Fr\'echet derivatives of $F$ with respect to $\eta$, and $c$, $b$ of all orders $\geq 3$ are zero everywhere by brutal force, and since $\c$ is a real analytic function, $F$ is a real analytic operator.

\subsection*{Bifurcation condition}

For any $T\geq 0$, $\k>0$, for any $c>0$, $b\in \mathbb{R}$ and $|b|$ sufficiently small, note that 
\begin{equation}\label{E:eta0'}
\eta_0(c;b,\k,T)=b(1-c)+O(b^2)
\end{equation}
makes a constant solution of \eqref{def:F}-\eqref{E:F=0} and, hence, \eqref{E:pFDCH}. It follows from the implicit function theorem that if non-constant solutions of \eqref{def:F}-\eqref{E:F=0} and, hence, \eqref{E:pFDCH} bifurcate from $\eta=\eta_0$ for some $c=c_0$ then, necessarily,
\[
L_0:=F_\eta(\eta_0,c_0;b,\k,T):H^{k+2}(\mathbb{T}) \to H^k(\mathbb{T}),
\]
where
\begin{equation}\label{def:L0}
L_0=\c(\k|\partial_z|;T) -c_0 -3 +3(1+\eta_0)^{1/2}+\frac{5}{12}\k^2\eta_0\partial_z^2,
\end{equation}
is not an isomorphism. Here $\eta_0$ depends on $c_0$. But we suppress it for simplicity of notation. A straightforward calculation reveals that $L_0e^{inz}=0$, $n\in\mathbb{Z}$, if and only if
\begin{equation}\label{E:bifurcation}
c_0=\c(n\k;T)-3+3(1+\eta_0)^{1/2}-\frac{5}{12}\k^2 n^2\eta_0.
\end{equation}
For $b=0$ and, hence, $\eta_0=0$ by \eqref{E:eta0'}, it simplifies to $c_0=\c(n\k;T)$. Without loss of generality, we restrict the attention to $n=1$. For $|b|$ sufficiently small, \eqref{E:bifurcation} and \eqref{E:eta0'} become
\begin{align}
c_0(b,\k,T)=&\c(\k;T)+b\Big(\frac32-\frac{5}{12}\k^2\Big)(1-c_{\rm ww}(\k;T))+O(b^2)\label{E:c0}
\intertext{and}
\eta_0(b,\k,T)=&b(1-c_{\rm ww}(\k;T))+O(b^2).\label{E:eta0}
\end{align}

For $T=0$, since $\c(\k;0)>\c(n\k;0)$ for $n=2,3,\dots$ everywhere in $\mathbb{R}$ (see Figure~\ref{fig:cT}a), it is straightforward to verify that for any $\k>0$, $b\in\mathbb{R}$ and $|b|$ sufficiently small, the kernel of $L_0:H^{k+2}(\mathbb{T})\to H^k(\mathbb{T})$ is two dimensional and spanned by $e^{\pm iz}$. Moreover, the co-kernel of $L_0$ is two dimensional. Therefore, $L_0$ is a Fredholm operator of index zero.

Similarly, for $T\geq1/3$, since $\c(\k;T)<\c(n\k;T)$ for $n=2,3,\dots$ everywhere in $\mathbb{R}$ (see Figure~\ref{fig:cT}b), for any $\k>0$, $b\in\mathbb{R}$ and $|b|$ sufficiently small, $L_0:H^{k+2}(\mathbb{T})\to H^k(\mathbb{T})$ is a Fredholm operator of index zero, whose kernel is two dimensional and spanned by $e^{\pm iz}$. 

For $0<T<1/3$, on the other hand, for any integer $n\geq 2$, it is possible to find some $\k$ such that $\c(\k;T)=\c(n\k;T)$ (see Figure~\ref{fig:cT}c). If 
\begin{equation}\label{A:resonance}
\c(\k;T)\neq\c(n\k;T)\qquad\text{for any $n=2,3,\dots$}
\end{equation} 
then $L_0:H^{k+2}(\mathbb{T})\to H^k(\mathbb{T})$ is likewise a Fredholm operator of index zero, whose kernel is two dimensional and spanned by $e^{\pm iz}$.  But if $\c(\k;T)=\c(n\k;T)$ for some integer $n\geq 2$, resulting in the resonance of the fundamental mode and the $n$-th harmonic, then the kernel is four dimensional.  One may follow along the same line as the argument in \cite{Jones}, for instance, to construct a periodic traveling wave. But we do not pursue this here. 

\subsection*{Lyapunov-Schmidt procedure}

For any $T\geq0$, $\k>0$ satisfying \eqref{A:resonance}, $b\in\mathbb{R}$ and $|b|$ sufficiently small, we employ a Lyapunov-Schmidt procedure to construct non-constant solutions of \eqref{def:F}-\eqref{E:F=0} and, hence, \eqref{E:pFDCH} bifurcating from $\eta=\eta_0$ and $c=c_0$, where $\eta_0$ and $c_0$ are in \eqref{E:eta0} and \eqref{E:c0}. Throughout the proof, $T$, $\k$, and $b$ are fixed and suppressed for simplicity of notation. 

Recall that $F(\eta_0,c_0)=0$, where $F$ is in \eqref{def:F}, and $L_0e^{\pm iz}=0$, where $L_0$ is in \eqref{def:L0}. We write that
\begin{equation}\label{def:eta2}
\eta(z)=\eta_0+\frac12(ae^{iz}+a^*e^{-iz})+\eta_r(z)\quad\text{and}\quad c=c_0+c_r,
\end{equation}
and we require that $a\in \mathbb{C}$, $\eta_r\in H^{k+2}(\mathbb{T})$ be even and  
\begin{equation}\label{E:v-condition}
\l \eta_r, e^{\pm iz}\r_{L^2(\mathbb{T})}=0,
\end{equation}
and $c_r\in \mathbb{R}$. Substituting \eqref{def:eta2} into \eqref{def:F}-\eqref{E:F=0}, we use $F(\eta_0,c_0)=0$, $L_0e^{\pm iz}=0$, and we make an explicit calculation to arrive at
\begin{align}
L_0\eta_r=&(3(1+\eta_0)^{1/2}+c_r)\Big(\frac12(ae^{iz}+a^*e^{-iz})+\eta_r\Big)\notag \\
&-2\Big(1+\eta_0+\frac12(ae^{iz}+a^*e^{-iz})+\eta_r\Big)^{3/2} \notag\\
&-\frac{13}{48}\k^2 \Big(\frac{i}{2}(ae^{iz}-a^*e^{-iz})+\eta_r'\Big)^2 \notag\\
&-\frac{5}{12}\k^2\Big(\frac12(ae^{iz}+a^*e^{-iz})+\eta_r\Big)\Big(-\frac12(ae^{iz}+a^*e^{-iz})+\eta_r''\Big)\notag\\
=:&g(\eta_r;a,a^*,c_r). \label{def:g}
\end{align}
Here and elsewhere, the prime means ordinary differentiation. Note that 
\[
g:H^{k+2}(\mathbb{T})\times\C\times \C\times \R\to H^k(\mathbb{T}).
\] 
Recall that $F$ is a real analytic operator. Hence $g$ depends analytically on its arguments. Clearly, $g(0;0,0,c_r)=0$ for all $c_r\in\mathbb{R}$. 

Let $\Pi:L^2(\mathbb{T})\to \ker L_0$ denote the spectral projection, defined as 
\[ 
\Pi f(z)=\widehat{f}(1)e^{iz}+\widehat{f}(-1)e^{-iz}.
\]
Since $\Pi\eta_r=0$ by \eqref{E:v-condition}, we may rewrite \eqref{def:g} as
\begin{equation}\label{E:LS1}
L_0\eta_r=(1-\Pi)g(\eta_r;a,a^*,c_r)\quad \text{and} \quad 0=\Pi g(\eta_r;a,a^*,c_r).
\end{equation}
Moreover, for any $T\geq0$, and $\k>0$ satisfying \eqref{A:resonance}, note that $L_0$ is invertible on $(1-\Pi)H^k(\mathbb{T})$. Specifically,
\[
{L_0}^{-1}f(z)=\sum_{n\neq \pm 1}\frac{\widehat{f}(n)}{\c(\k n;T)-\c(\k;T)+\frac{5}{12}\k^2 \eta_0(1- n^2)}e^{inz}.
\]
Hence we may rewrite \eqref{E:LS1} as
\begin{equation}\label{E:LS2}
\eta_r=L_0^{-1}(I-\Pi)g(\eta_r;a,a^*,c_r)\quad \text{and}\quad 0=\Pi g(\eta_r;a,a^*,c_r).
\end{equation}
Note that $L_0^{-1}: (1-\Pi)H^k(\mathbb{T}) \to H^k(\mathbb{T})$ is bounded. We claim that
\[
L_0^{-1}: (1-\Pi)H^k(\mathbb{T}) \to H^{k+2}(\mathbb{T})
\]
is bounded. As a matter of fact, 
\[
\left|\frac{n^2\widehat{f}(n)}{\c(\k n;T)-\c(\k;T)+\frac{5}{12}\k^2 \eta_0(1- n^2)}\right|\leq C|\widehat{f}(n)|
\]
for some constant $C>0$ for $n\in\mathbb{Z}$ and $|n|$ sufficiently large. Therefore, for any $a, a^*\in\mathbb{C}$ and $c_r\in\mathbb{R}$, 
\[
L_0^{-1}(1-\Pi)g:H^{k+2}(\mathbb{T})\to H^{k+2}(\mathbb{T})
\]
is bounded. Note that it depends analytically on its argument. Since $g(0;0,0,c_r)=0$ for any $c_r\in\mathbb{R}$, it follows from the implicit function theorem that a unique solution 
\[
\eta_2=\eta_r(a,a^*,c_r)
\] 
exists to the former equation of \eqref{E:LS2} near $\eta_r=0$ for $a\in\mathbb{C}$ and $|a|$ sufficiently small for any $c_r\in\mathbb{R}$. Note that $\eta_2$ depends analytically on its arguments and it satisfies \eqref{E:v-condition} for $|a|$ sufficiently small for any $c_r\in\mathbb{R}$. The uniqueness implies
\begin{equation}\label{Vdef}
\eta_2(0,0,c_r)=0 \qquad\text{for any $c_r\in \mathbb{R}$}.
\end{equation}
Moreover, since \eqref{def:F}-\eqref{E:F=0} and, hence, \eqref{E:LS2} are invariant under \eqref{E:invariance} for any $z_0\in\mathbb{R}$, it follows that 
\begin{align}
\eta_2(a,a^*,c_r)(z+z_0)=\eta_2(ae^{iz_0},a^*e^{-iz_0},c_r)\quad\text{and}\quad 
\eta_2(a,a^*,c_r)(-z)=\eta_2(a,a^*,c_r)(z) \label{Vprop}
\end{align}
for any $z_0\in\mathbb{R}$ for any $a\in\mathbb{C}$ and $|a|$ sufficiently small, and $c_r\in\mathbb{R}$. 

To proceed, we rewrite the latter equation in \eqref{E:LS2} as 
\[
\Pi g(\eta_2(a,a^*,c_r);a,a^*,c_r)=0
\]
for $a\in\mathbb{C}$ and $|a|$ sufficiently small for $c_r\in\mathbb{R}$. This is solvable, provided that
\begin{equation}\label{def:pi}
\pi_{\pm}(a,a^*,c_r):=\left\l g(\eta_2(a,a^*,c_r);a,a^*,c_r),ae^{iz}\pm a^*e^{-iz}\right\r_{L^2(\mathbb{T})}=0.
\end{equation}
We use \eqref{Vprop}, where $z_0=-2\arg (a)$, and \eqref{def:pi} to show that 
\[
\pi_- (a^*,a,c_r)=\pi_- (a,a^*,c_r)=-\pi_- (a^*,a,c_r).
\]
Hence $\pi_- (a,a^*,c_r)=0$ holds for any $a\in\mathbb{C}$ and $|a|$ sufficiently small for any $c_r\in\mathbb{R}$. Moreover, we use \eqref{Vprop}, where $z_0=-\arg (a)$, and \eqref{def:pi} to show that 
\[
\pi_+ (a,a^*,c_r)=\pi_+ (|a|,|a|,c_r).
\] 
Hence it suffices to solve $\pi_+ (a,a,c_r)=0$ for any $a, c_r\in \mathbb{R}$ and $|a|$ sufficiently small. 

Substituting \eqref{def:g} into \eqref{def:pi}, where $\eta_r=\eta_2(a,a,c_r)$, we make an explicit calculation to arrive at 
\[
\pi_+ (a,a,c_r)=a^2(\pi c_r+\pi_r(a,c_r)),
\]
where 
\begin{align*}
\pi_r(a,c_r)=&-2a^{-1} \l(1+\eta_0+a\cos z+\eta_2(a,a,c_r)(z))^{3/2},\cos z\r\\
&-\frac{5}{12}\k^2(\l \eta_2''(a,a,c_r)(z) - \eta_2(a,a,c_r)(z),\cos^2 z\r
-a^{-1}\l \eta_2\eta_2''(a,a,c_r)(z), \cos z\r)\\
& -\frac{13}{48}\k^2 a^{-1}( \l \eta_2'(a,a,c_r)(z)^2,  \cos z\r-\l\eta_2'(a,a,c_r)(z),\sin 2z\r),
\end{align*}
and $\l\cdot\,,\cdot\r$ means the $L^2(\mathbb{T})$ inner product. We merely pause to remark that $\pi_r$ is well defined. As a matter of fact, $a^{-1}\eta_2$ is not singular for $a\in\mathbb{R}$ and $|a|$ sufficiently small by \eqref{Vdef}. Clearly, $\pi_r$ and, hence, $\pi_{\pm}$ depend analytically on its arguments. Since $\pi_r(0,0)=\partial \pi_r/\partial c_r(0,0)=0$ by \eqref{Vdef}, it follows from the implicit function theorem that a unique solution 
\[
c_r=c_1(a)
\] 
exists to $\pi_+(a,a,c_r)=0$ and, hence, the latter equation of \eqref{E:LS2} near $c_r=0$ for $a\in\mathbb{R}$ and $|a|$ sufficiently small. Clearly, $c_1$ depends analytically on $a$. 

To recapitulate,
\[ 
\eta_r=\eta_2(a,a,c_1(a))\quad\text{and}\quad c_r=c_1(a)
\]
uniquely solve \eqref{E:LS2} for $a\in\mathbb{R}$ and $|a|$ sufficiently small, and by virtue of \eqref{def:eta2}, 
\begin{equation}\label{E:soln}
\eta(a)(z)=\eta_0+a\cos z+\eta_2(a,a,c_1(a))(z)\quad\text{and}\quad  c(a)=c_0+c_1(a)
\end{equation}
uniquely solve \eqref{def:F}-\eqref{E:F=0} and, hence, \eqref{E:pFDCH} for $a\in\mathbb{R}$ and $|a|$ sufficiently small. Note that $\eta$ is $2\pi$ periodic and even in $z$. Moreover, $\eta \in H^\infty(\mathbb{T})$. 

For $a,b\in\mathbb{R}$ and $|a|,|b|$ sufficiently small, we write that 
\begin{align}
\eta(a;b,\k,T)(z):=&\eta_0(b,\k,T)+ a \cos z+a^2\eta_2(z)+a^3\eta_3(z)+\cdots\label{E:eta(a,b)}
\intertext{and} 
c(a;b,\k,T):=&c_0(b,\k,T)+ac_1+a^2c_2+\cdots,\label{E:c(a,b)}
\end{align}
where $\eta_2,\eta_3, \dots$ are $2\pi$ periodic, even, and smooth functions of $z$, and $c_1, c_2, \dots \in \mathbb{R}$. 

We claim that $c_1=0$. As a matter of fact, note that \eqref{E:pFDCH} and, hence, \eqref{def:F}-\eqref{E:F=0} remain invariant under $z\mapsto z+\pi$ by \eqref{E:invariance}. Since $\partial \eta/\partial a(0)(z)=\cos z$, however, $\eta(z)\neq\eta(z+\pi)$ must hold. Thus $\partial c/\partial a(0)=0$. This proves the claim. If $\l\eta_{j-1},\eta_j\r_{L^2(\mathbb{T})}=0$ for any integer $j\geq 1$, in addition, then $c_{2j-1}=0$ for any integer $j\geq1$. Hence $c$ is even in $a$.

Substituting \eqref{E:eta(a,b)} and \eqref{E:c(a,b)} into \eqref{E:pFDCH}, we may calculate the small amplitude expansion. The proof is very similar to that in \cite{HP1}, for instance. Hence we omit the details.

Below we summarize the conclusion.

\begin{lemma}[Existence of sufficiently small and periodic traveling waves]\label{lem:existence} 
For any $T\geq0$, $\k>0$ satisfying \eqref{A:resonance}, $b\in\mathbb{R}$ and $|b|$ sufficiently small, a one parameter family of solutions of \eqref{E:pFDCH} exists, denoted $\eta(a;b,\k,T)$ and $c(a;b,\k,T)$, for $a\in\mathbb{R}$ and $|a|$ sufficiently small; $\eta\in H^\infty(\mathbb{T})$ and it is even in $z$; $\eta$ and $c$ depend analytically on $a$, and $b$, $\k$. Moreover,
\begin{align}
\eta(a;b,\k,T)(z)=&b(1-c_{\rm ww}(\k;T))+a\cos z+a^2(h_0+h_2\cos 2z)+O(a(a+b)^2),\hspace*{-20pt}\label{E:n(k,a,b)} \\
c(a;b,\k,T)=&c_{\rm ww}(\k;T)+b\Big(\frac32-\frac{5}{12}\k^2\Big)(1-c_{\rm ww}(\k;T))+a^2 c_2+O(a(a+b)^2) \hspace*{-20pt}\label{E:c(k,a,b)}
\end{align}
as $a,b \to 0$, where
\begin{gather} 
h_0=\Big(\frac38-\frac{7}{96}\k^2\Big)\frac{1}{c_{\rm ww}(\k;T)-1},\qquad
h_2 =\Big(\frac38-\frac{11}{32}\k^2\Big)\frac{1}{c_{\rm ww}(\k;T)-c_{\rm ww}(2\k;T)}, \label{def:A0A2}
\intertext{and}
c_2=\Big(\frac32-\frac{5}{12}\k^2\Big)h_0+\Big(\frac34-\frac12\k^2\Big)h_2-\frac{3}{32}. \label{def:c2}
\end{gather}
\end{lemma}

\section{Modulational instability index}\label{sec:MI}

For $T\geq0$, $\k>0$ satisfying \eqref{A:resonance}, $a, b\in\mathbb{R}$ and $|a|, |b|$ sufficiently small, let $\eta=\eta(a;b,\k,T)$ and $c=c(a;b,\k,T)$,  denote a sufficiently small and $2\pi/\k$ periodic traveling wave of \eqref{E:FDCH} and \eqref{def:cT1}, whose existence follows from the previous section. We address its modulational stability and instability. 

Linearizing \eqref{E:FDCH} about $\eta$ in the coordinate frame moving at the speed $c$, we arrive at
\[
\zeta_t+\k\partial_z\Big(\c(\k |\partial_z|;T) -c-3+3(1+\eta)^{1/2}+ \k^2\Big(\frac{5}{12}(\eta \partial_z^2+\eta_{zz})+\frac{13}{24}  \eta_z \partial_z\Big)\Big)\zeta=0,
\]
where $\c(\k|\partial_z|;T)$ is in \eqref{def:cT1}. Seeking a solution of the form $\zeta(z,t)=e^{\lambda \k t}\zeta(z)$, $\lambda\in\mathbb{C}$, we arrive at
\begin{align}\label{E:eigen}
\lambda \zeta&=\partial_z\Big(-\c(\k |\partial_z|;T)+c+3-3(1+\eta)^{1/2}-  \k^2\Big(\frac{5}{12}(\eta \partial_z^2+\eta_{zz})+\frac{13}{24} \eta_z \partial_z\Big)\Big)\zeta \\
&=:\mathcal{L}(a;b,\k,T)\zeta. \notag
\end{align}
We say that $\eta$ is {\em spectrally unstable} to square integrable perturbation if the $L^2(\mathbb{R})$ spectrum of $\L$ intersects the open right-half plane of $\mathbb{C}$, and it is spectrally stable otherwise. Note that $\eta$ is $2\pi$ periodic in $z$, but $\zeta$ needs not. Note that the spectrum of $\L$ is symmetric with respect to the reflections in the real and imaginary axes. Hence $\eta$ is spectrally unstable if and only if the spectrum of $\L$ is {\em not} contained in the imaginary axis. 

It is well known (see \cite{BHJ}, for instance, and references therein) that the $L^2(\mathbb{R})$ spectrum of $\L$ contains no eigenvalues. Rather, it consists of the essential spectrum. Moreover, a nontrivial solution of \eqref{E:eigen} does not belong to $L^p(\mathbb{R})$ for any $p\in [1,\infty)$. Rather, if $\zeta\in L^\infty(\mathbb{R})$ solves \eqref{E:eigen} then, necessarily,
\[
\zeta(z)=e^{i\xi z}\phi(z),\qquad\text{where}\quad \phi(z+2\pi)=\phi(z),
\]
for some $\xi$ in the range $(-1/2,1/2]$. It follows from Floquet theory (see \cite{BHJ}, for instance, and references therein) that $\lambda$ belongs to the $L^2(\mathbb{R})$ spectrum of $\mathcal{L}$ if and only if 
\begin{equation}\label{def:Lx}
\lambda\phi=e^{-i\xi z}\mathcal{L}(a;b,\k,T)e^{i\xi z}\phi=:\mathcal{L}(\xi)(a;b,\k,T)\phi
\end{equation}
for some $\xi\in(-1/2,1/2]$ and $\phi\in L^2(\mathbb{T})$. Thus
\[
\text{spec}_{L^2(\mathbb{R})}(\mathcal{L}(a;b,\k,T))=
\bigcup_{\xi\in(-1/2,1/2]}\text{spec}_{L^2(\mathbb{T})}(\mathcal{L}(\xi)(a;b,\k,T)).
\]
Note that for any $\xi\in(-1/2,1/2]$, the $L^2(\mathbb{T})$ spectrum of $\mathcal{L}(\xi)$ comprises of eigenvalues of finite multiplicities. Thus the essential spectrum of $\L$ may be characterized as a one parameter family of point spectra of $\L(\xi)$ for $\xi\in(-1/2,1/2]$. Note that 
\[
\text{spec}_{L^2(\mathbb{T})}(\L(\xi))=(\text{spec}_{L^2(\mathbb{T})}(\L(-\xi)))^*.
\]
Hence it suffices to take $\xi\in [0,1/2]$.

Note that $\xi=0$ corresponds to the same period perturbations as $\eta$. Moreover, $\xi>0$ and small corresponds to long wavelength perturbations, whose effects are to slowly vary the period and other wave characteristics. They furnish the spectral information of $\L$ in the vicinity of the origin in $\mathbb{C}$; see \cite{BHJ}, for instance, for details. Therefore, we say that $\eta$ is {\em modulationally unstable} if the $L^2(\mathbb{T})$ spectra of $\L(\xi)$ are not contained in the imaginary axis near the origin for $\xi>0$ and small, and it is modulationally stable otherwise. 

For an arbitrary $\xi$, one must in general study \eqref{def:Lx} by means of numerical computation. But, for $\xi>0$ and small for $\lambda$ in the vicinity of the origin in $\mathbb{C}$, we may take a spectral perturbation approach in \cite{HJ2, HJ3, HP1}, for instance, to address it analytically. 

\subsubsection*{Notation}
Throughout the section, $T\geq0$ and $\k>0$ satisfying \eqref{A:resonance} are suppressed for simplicity of notation, unless specified otherwise. We assume that $b=0$. For nonzero $b$, one may explore in like manner. But the calculation becomes lengthy and tedious. Hence we do not discuss the details. We use
\begin{equation}\label{def:Lxia}
\L(\xi,a)=\L(\xi)(a;0,\k,T).
\end{equation}

\subsection*{Spectra of $\L(\xi,0)$ and $\L(0,a)$}

For $a=0$ --- namely, the rest state --- a straightforward calculation reveals that 
\begin{equation} \label{eq:a=0FDCH}
\L(\xi,0)e^{inz}=i\omega(n+\xi)e^{inz}
\quad \text{for $n\in\mathbb{Z}$ and $\xi\in[0,1/2]$},
\end{equation}
where 
\begin{equation}\label{def:w-FDCH}
\omega(n+\xi)=(\xi+n)(c_{\rm ww}(\k;T)-\c(\k(n+\xi);T)).
\end{equation}
For $\xi=0$, 
\[
\omega(1)=\omega(-1)=\omega(0)=0,
\]
and $\omega(n)\neq 0$ otherwise. Hence zero is an $L^2(\mathbb{T})$ eigenvalue of $\mathcal{L}(0,0)$ with  multiplicity three. Moreover,
\begin{equation}\label{eq:eigen0-FDCH}
\cos z, \qquad \sin z,\quad\text{and}\quad 1
\end{equation}
are the associated eigenfunctions, real valued and orthogonal to each other. For $\xi>0$ sufficiently small, 
\[
i\omega(\pm1+\xi)\quad\text{and}\quad i\omega(\xi)
\]
are the $L^2(\mathbb{T})$ eigenvalues of $\L(\xi,0)$ in the vicinity of the origin in $\mathbb{C}$, and \eqref{eq:eigen0-FDCH} are the associated eigenfunctions. 

For $a\in\mathbb{R}$ and $|a|$ sufficiently small for $\xi=0$, zero is an $L^2(\mathbb{T})$ eigenvalue of $\mathcal{L}(0,a)$ with algebraic multiplicity three and geometric multiplicity two, and 
\begin{equation}\label{def:p123}
\begin{aligned}
\phi_1(z) &:=\Big(\eta_a-\frac{c_a}{c_b}\eta_b\Big)(a;0,\k,T)(z)=\cos z+ap_1+2a h_2\cos 2z+O(a^2),\\
\phi_2(z) &:=-\frac{1}{a}\eta_z(a;0,\k,T)(z)= \sin z+2a h_2\sin 2z+O(a^2),\\
\phi_3(z) &:=\frac{1}{1-\c(\k;T)}\eta_b(\k,a,0;T)(z)= 1+O(a^2)
\end{aligned}
\end{equation}
are the associated eigenfunctions, where 
\begin{equation}\label{def:p1}
p_1 = 2h_0-\frac{24c_2}{18-5 \k^2} = \frac{1}{18-5 \k^2}\Big( \frac94 -\frac{3}{16}\frac{(3-2 \k^2)(12-11 \k^2)}{c_{\rm ww}(\k;T)-c_{\rm ww}(2\k;T)}\Big)
\end{equation} 
and $h_2$ is defined in \eqref{def:A0A2}. The proof is nearly identical to that in \cite{HJ2}, for instance. Hence we omit the details. For $a=0$, note that \eqref{def:p123} becomes \eqref{eq:eigen0-FDCH}.

\subsection*{Spectral perturbation calculation}

Recall that for $\xi>0$ and sufficiently small for $a=0$, the $L^2(\mathbb{T})$ spectrum of $\L(\xi,0)$ contains three purely imaginary eigenvalues $i\omega(\pm1+\xi)$ and $i\omega(\xi)$ in the vicinity of the origin in $\mathbb{C}$, and \eqref{eq:eigen0-FDCH} spans the associated eigenspace, which does not depend on $\xi$. For $\xi=0$ for $a\in\mathbb{R}$ and $|a|$ sufficiently small, the spectrum of $\L(0,a)$ contains three eigenvalues at the origin, and \eqref{def:p123} spans the associated eigenspace, which depends analytically on $a$. 

For $\xi>0$, $a\in\mathbb{R}$ and $\xi, |a|$ sufficiently small, it follows from perturbation theory (see \cite{K}, for instance, for details) that the $L^2(\mathbb{T})$ spectrum of $\L(\xi,a)$ contains three eigenvalues in the vicinity of the origin in $\mathbb{C}$, and \eqref{def:p123} spans the associated eigenspace. Let
\begin{equation}\label{def:BI}
\mathbf{L}(\xi,a)=\left( \frac{\l \L(\xi,a)\phi_j, \phi_k\r}{\l \phi_j, \phi_j\r}\right)_{j,k=1,2,3}
\quad\text{and}\quad
\mathbf{I}(a)=\left( \frac{\l \phi_j, \phi_k\r}{\l \phi_j, \phi_j\r}\right)_{j,k=1,2,3},
\end{equation}
where $\phi_1,\phi_2,\phi_3$ are in \eqref{def:p123}. Throughout the subsection, $\langle\cdot\,,\cdot\rangle$ means the $L^2(\mathbb{T})$ inner product. Note that $\mathbf{L}$ represents the action of $\mathcal{L}$ on the eigenspace, spanned by $\phi_1,\phi_2,\phi_3$, and $\mathbf{I}$ is the projection of the identity onto the eigenspace. It follows from perturbation theory (see \cite{K}, for instance for details) that for $\xi>0$, $a\in\mathbb{R}$ and $\xi, |a|$ sufficiently small, the eigenvalues of $\L(\xi,a)$ agree in location and multiplicity with the roots of $\det(\mathbf{L}-\lambda\mathbf{I})$ up to terms of order $a$.

For $a\in\mathbb{R}$ and $|a|$ sufficiently small, a Baker-Campbell-Hausdorff expansion reveals that
\[
\L(\xi,a)=\L(0,a)+i\xi [\L(0,a),z]-\frac12\xi^2[[\L(0,a),z],z]+O(\xi^3)
\]
as $\xi\to0$, where $[\cdot\,,\cdot]$ means the commutator. We merely pause to remark that $[\L,z]$ and $[[\L,z],z]$ are well defined in the periodic setting even though $z$ is not. We use \eqref{E:eigen}, \eqref{def:Lx} and \eqref{E:n(k,a,b)}, \eqref{E:c(k,a,b)} to write
\begin{align}
\L(\xi,a)=&\mathcal{M}
-a\partial_z\Big(\frac32\cos z+\k^2\Big(\frac{5}{12}\cos z(\partial_z^2-1)-\frac{13}{24}\sin z\partial_z\Big)\Big) \label{def:L-FDCH}\\
&-i\xi a\Big(\frac32\cos z  + \k^2\Big(\frac{5}{12}(2\partial_z\cos z\partial_z+\cos z(\partial_z^2-1))
-\frac{13}{24} (\sin z\partial_z+\partial_z \sin z)\Big)\Big)\notag \\&+O(\xi^3+\xi^2 a+a^2)\notag
\end{align}
as $\xi, a\to 0$, where 
\[
\mathcal{M}=\L(0,0)+i\xi [\L(0,0),z]-\frac12\xi^2[[\L(0,0),z],z]
\]
agrees with $\L(\xi,0)$ up to terms of order $\xi^2$ as $\xi\to 0$. We may then resort to \eqref{eq:a=0FDCH}, \eqref{def:w-FDCH}, and we make an explicit calculation to find that 
\begin{align*}
\L(\xi,0) e^{ inz}=&in(c_{\rm ww}(\k;T)-c_{\rm ww}(n\k;T))e^{ inz}\\
&+i\xi(c_{\rm ww}(\k;T)-c_{\rm ww}(n\k;T)-\k c'_{\rm ww}(n\k;T))e^{ inz} \\
& -\frac12\xi^2(2\k c'_{\rm ww}(n\k;T)+\k^2 c''_{\rm ww}(n\k;T))e^{ inz}+O(\xi^3)
\end{align*}
as $\xi\to 0$.
Therefore, for $T\geq0$ but $T\neq1/3$, for $\k>0$ satisfying \eqref{A:resonance}, $\M 1=i\xi(c_{\rm ww}(\k;T)-1)$, 
\begin{align*}
\M\left\{ \begin{matrix} \cos z\\ \sin z\end{matrix}\right\}=&
-i\xi \k \c'(\k;T)\left\{ \begin{matrix} \cos z\\ \sin z\end{matrix}\right\} 
\pm \frac12\xi^2(2\k \c'(\k;T)+\k^2\c''(\k;T))\left\{ \begin{matrix} \sin z\\ \cos z\end{matrix}\right\},
\end{align*} 
and
\begin{align*}
\M\left\{ \begin{matrix} \cos 2z\\ \sin 2z\end{matrix}\right\}=&
\mp 2(\c(\k;T)-\c(2\k;T))\left\{ \begin{matrix} \sin 2z\\ \cos 2z\end{matrix}\right\}\\
&+i\xi(\c(\k;T)-\c(2\k;T)-2\k c'_{\rm ww}(2\k;T))\left\{ \begin{matrix} \cos 2z\\ \sin 2z\end{matrix}\right\} \\
&\pm \frac12\xi^2(2\k c'_{\rm ww}(2\k;T)+\k^2 c''_{\rm ww}(2\k;T))\left\{ \begin{matrix} \sin 2z\\ \cos 2z\end{matrix}\right\}.
\end{align*} 
For $T=1/3$, one must calculate higher order terms in the expansion of $\L(\xi,0)$. We do not pursue this here. 

We use \eqref{def:L-FDCH}, \eqref{def:p123} and the above formula for $\mathcal{M}$, and we make a lengthy but explicit calculation to find that
\begin{align*}
\L\phi_1=&-i\xi \k \c'(\k;T)\cos z \\ 
&-i\xi a\Big(\frac34-\frac{7}{48}\k^2-p_1(\c(\k;T)-1)\Big)\\
&+i\xi a\Big(-\frac34+\frac{33}{16}\k^2+2h_2(\c(\k;T)-\c(2\k;T)-2\k c'_{\rm ww}(2\k;T))\Big)\cos 2z\notag \\
&+\frac12\xi^2(2\k  \c'(\k;T)+\k^2\c''(\k;T))\sin z+O(\xi^3+\xi^2a+a^2), \\
\L\phi_2 =&-i\xi \k \c'(\k;T)\sin z \\
&+i\xi a\Big(-\frac34+\frac{33}{16}\k^2+2h_2(\c(\k;T)-\c(2\k;T)-2\k c'_{\rm ww}(2\k;T))\Big)\sin 2z\notag\\
&-\frac12\xi^2(2\k \c'(\k;T)+\k^2 \c''(\k;T))\cos z+O(\xi^3+\xi^2a+a^2), \\
\L\phi_3=&a\Big(3-\frac{5}{6}\k^2\Big)\sin z+i\xi (\c(\k;T)-1)-i\xi a\Big(\frac32-\frac{23}{24}\k^2\Big)\cos z+O(\xi^3+\xi^2a+a^2) 
\end{align*}
as $\xi, a\to 0$, where $p_1$ is in \eqref{def:p1} and $h_2$ is in \eqref{def:A0A2}.

To proceed, we take the $L^2(\mathbb{T})$ inner product of the above and \eqref{def:p123}, and we make a lengthy but explicit calculation to find that
\begin{align*}
\langle \L\phi_1,\phi_1\rangle &=\langle\L\phi_2,\phi_2 \rangle
= -\frac12i\xi \k \c'(\k;T)+O(\xi^3+\xi^2a+a^2),\\
\langle \L\phi_1,\phi_2\rangle &=-\langle \L\phi_2,\phi_1\rangle= \frac14\xi^2(2\k \c'(\k;T)+\k^2 \c''(\k;T))+O(\xi^3+\xi^2a+a^2),\\
\langle \L\phi_1,\phi_3\rangle &= 
i\xi a\Big(-\frac34+\frac{36}{48}\k^2+p_1 (\c(\k;T)-1)\Big)+O(\xi^3+\xi^2a+a^2),
\intertext{and}
\langle \L\phi_2,\phi_3\rangle &=0+O(\xi^3+\xi^2a+a^2),\\
\langle \L\phi_3,\phi_1\rangle &=i\xi a\Big(-\frac34+\frac{23}{48}\k^2+p_1 (\c(\k;T)-1) \Big)+O(\xi^3+\xi^2a+a^2),\\
\langle \L\phi_3,\phi_2\rangle &=a\Big(\frac34-\frac{5}{24}\k^2\Big)+O(\xi^3+\xi^2a+a^2),\\
\langle \L\phi_3,\phi_3\rangle &=i\xi( \c(\k;T)-1)+O(\xi^3+\xi^2a+a^2)
\end{align*}
as $\xi,a\to 0$,  where $p_1$ is in \eqref{def:p1}. Moreover, we take the $L^2(\mathbb{T})$ inner products of \eqref{def:p123}, and we make an explicit calculation to find that 
\begin{align*}
\langle \phi_1,\phi_1\rangle &= \langle \phi_2,\phi_2 \rangle=\frac{1}{2}+O(\xi^3+\xi^2a+a^2),\\
\langle \phi_1,\phi_2\rangle &, \langle \phi_2,\phi_3\rangle=0+O(\xi^3+\xi^2a+a^2),\\
\langle \phi_1,\phi_3\rangle &=ap_1+O(\xi^3+\xi^2a+a^2), \\
\langle \phi_3,\phi_3\rangle &=1+O(\xi^3+\xi^2a+a^2)
\end{align*}
as $\xi,a\to 0$, where $p_1$ is in \eqref{def:p1}. Together, \eqref{def:BI} becomes
\begin{align}\label{eq:B-FDCH}
\mathbf{L}(\xi,a) &=a\Big(\frac34-\frac{5}{24}\k^2\Big)\begin{pmatrix} 0&0&0\\0&0&0\\0&1&0\end{pmatrix}\\
&+i\xi \begin{pmatrix} -\k \c'(\k;T)&0&0\\ 0&-\k \c'(\k;T)&0\\ 0&0&\c(\k;T)-1\end{pmatrix} \notag\\
&+i\xi a\begin{pmatrix} 0&0&{\displaystyle -\tfrac34+\tfrac{7}{24} \k^2+2p_1 (\c(\k;T)-1)}\\ 0&0&0\\ 
{\displaystyle -\tfrac34+\tfrac{23}{48}\k^2+p_1 (\c(\k;T)-1)}&0&0\end{pmatrix}\notag\\
&+\xi^2(\k \c'(\k;T)+\tfrac12\k^2\c''(\k;T))\begin{pmatrix}0&1&0\\ -1&0&0\\ 0&0&0 \end{pmatrix}
+O(\xi^3+\xi^2a+a^2),\notag
\end{align}
and
\begin{equation}
\mathbf{I}(a)=\mathbf{I}+ap_1 \begin{pmatrix} 0&0&2\\0&0&0\\1&0&0\end{pmatrix}+O(a^2)\label{eq:I-FDCH}
\end{equation}
as $\xi,a\to 0$, where $p_1$ is in \eqref{def:p1} and $\mathbf{I}$ is the $3\times3$ identity matrix. 

\subsection*{Modulational instability index}

For any $T\geq0$ but $T\neq1/3$, $\k>0$ satisfying \eqref{A:resonance}, for $\xi>0$, $a\in\mathbb{R}$ and $\xi, |a|$ sufficiently small, we turn the attention to the roots of 
\begin{multline*}
\det(\mathbf{L}-\lambda\mathbf{I})(\xi,a;\k,T) \\
=:p_3(\xi,a;\k,T)\lambda^3+ip_2(\xi,a;\k,T)\lambda^2+p_1(\xi,a;\k,T)\lambda+ip_0(\xi,a;\k,T),
\end{multline*}
where $\mathbf{L}$ and $\mathbf{I}$ are in \eqref{eq:B-FDCH} and \eqref{eq:I-FDCH}.
Details are found in \cite{HJ2}, for instance. Hence we merely hit the main points.

Let 
\[
q(-i\xi\lambda)(\xi,a;\k,T)=(i\xi^3(q_3\lambda^3-q_2\lambda^2-q_1\lambda+q_0)(\xi,a;\k,T),
\]
where $p_j=\xi^{3-j}q_j$, $j=0,1,2,3$. Note that $q_0,q_1,\dots,q_3$ are real valued and depend analytically on $\xi,a$ and $\k$ for any $\xi>0$ and $|a|$ sufficiently small for any $\k>0$ satisfying \eqref{A:resonance}. Moreover, they are odd in $\xi$ and even in $a$. For any $T\geq0$ but $T\neq1/3$, $\k>0$ satisfying \eqref{A:resonance}, $a\in\mathbb{R}$ and $|a|$ sufficiently small, a periodic traveling wave $\eta(a;0,\k,T)$ and $c(a;0,\k,T)$ of \eqref{E:FDCH} and \eqref{def:cT1} is modulationally unstable, provided that $q$ possesses a pair of complex roots or, equivalently, 
\[
\Delta_0(\xi,a;\k,T):=(18q_3q_2q_1q_0+q_2^2q_1^2+4q_2^3q_0+4q_3q_1^3-27q_3^2q_0^2)(\xi,a;\k,T)<0
\]
for $\xi>0$ and small, and it is modullationally stable if $\Delta_0>0$. Note that $\Delta_0$ is even in $\xi$ and $a$. Hence we write that 
\[
\Delta_0(\xi,a;\k,T)=:\Delta_0(\xi,0;\k,T)+a^2\Delta(\k;T)+O(a^2(\xi^2+a^2))
\]
as $a\to 0$ for $\xi>0$ and small. We then use \eqref{eq:B-FDCH} and \eqref{eq:I-FDCH}, and we make a Mathematica calculation to show that
\[
\Delta_0(\xi,0;\k,T)=\k^2 i_1^2(\k;T)\Big(\xi i_2^2+\frac14\xi^3\k^2 i_1^2\Big)^2(\k;T)>0
\]
as $\xi\to0$. Therefore, if $\Delta<0$ then $\Delta_0<0$ for $\xi>0$ and sufficiently small, depending on $a\in\mathbb{R}$ and $|a|$ sufficiently small, implying modulational instability, whereas if $\Delta>0$ then $\Delta_0>0$ 
for $\xi>0$, $a\in\mathbb{R}$ and $\xi, |a|$ sufficiently small, implying modulational stability. We use \eqref{eq:B-FDCH} and \eqref{eq:I-FDCH}, and we make a Mathematica calculation to find $\Delta$ explicitly. 

Below we summarize the conclusion.

\begin{theorem}[Modulational instability index]\label{thm:FDCH} 
 For any $T \geq 0$ but $T\neq 1/3$, for any $\k>0$ satisfying \eqref{A:resonance}, a sufficiently small and $2\pi/\k$ periodic traveling wave of \eqref{E:FDCH} and \eqref{def:cT1} is modulationally unstable, provided that 
\begin{equation} \label{def:index}
\Delta(\k;T):=\frac{i_1i_2 }{i_3}i_4(\k;T)<0,
\end{equation}
where
\begin{subequations}\label{def:i1234}
\begin{align}
i_1(\k;T)=&(\k c_{\rm ww}(\k;T))'',  \\
i_2(\k;T)=&(\k c_{\rm ww}(\k;T))'-1,\\
i_3(\k;T)=&c_{\rm ww}(\k;T)-c_{\rm ww}(2\k;T), \\
i_4(\k;T)=&\Big(3i_2-i_2i_3+6 i_3
-\frac{1}{12}\k^2(57i_2+34 i_3)+\frac{1}{108}\k^4(198 i_2+35i_3)\Big)(\k,T), \hspace*{-30pt}\label{def:i4}
\end{align}
\end{subequations}
and $\c(\k;T)$ is in \eqref{def:cT1}. It is modulationally stable if $\Delta(\k;T)>0$.
\end{theorem}

Theorem~\ref{thm:FDCH} elucidates four resonance mechanisms which contribute to the sign change in $\Delta$ and, ultimately, the change in the modulational stability and instability in \eqref{E:FDCH} and \eqref{def:cT1}. Note that 
\[
\text{$\c(\k;T)=$the phase speed\quad and\quad $(\k\c(\k;T))'=$the group speed}
\]
in the linear theory. Specifically,  
\begin{itemize}
\item[(R1)] $i_1(\k;T)=0$ at some $\k$; the group speed achieves an extremum at the wave number $\k$; 
\item[(R2)] $i_2(\k;T)=0$ at some $\k$; the group speed at the wave number $\k$ coincides with the phase speed in the limit as $\k\to0$, resulting in the ``resonance of short and long waves;"
\item[(R3)] $i_3(\k;T)=0$ at some $\k$; the phase speeds of the fundamental mode and the second harmonic coincide at the wave number $\k$, resulting in the ``second harmonic resonance;"
\item[(R4)] $i_4(\k;T)=0$ at some $\k$.
\end{itemize}
Resonances (R1), (R2), (R3) are determined from the dispersion relation. For instance, $i_1, i_2, i_3$ appear in an index formula for the Whitham equation; see \cite{HJ3} for details. Resonance (R4), on the other hand, results from the resonance of the dispersion and nonlinear effects, and it depends on the nonlinearity of the equation. 

\section{Results}\label{sec:results}

\subsection*{For $\mathbf{T=0}$}

Since $(\k\c(\k;0))'<1$ for any $\k>0$ and decreases monotonically over the interval $(0,\infty)$ by brutal force, $i_1(\k;0)<0$ and $i_2(\k;0)<0$ for any $\k>0$. Since $\c(\k;0)>0$ for any $\k>0$ and decreases monotonically over the interval $(0,\infty)$ (see Figure~\ref{fig:cT}a), $i_3(\k;0)>0$ for any $\k>0$. 

\begin{figure}[h] 
~~\includegraphics[scale=0.7]{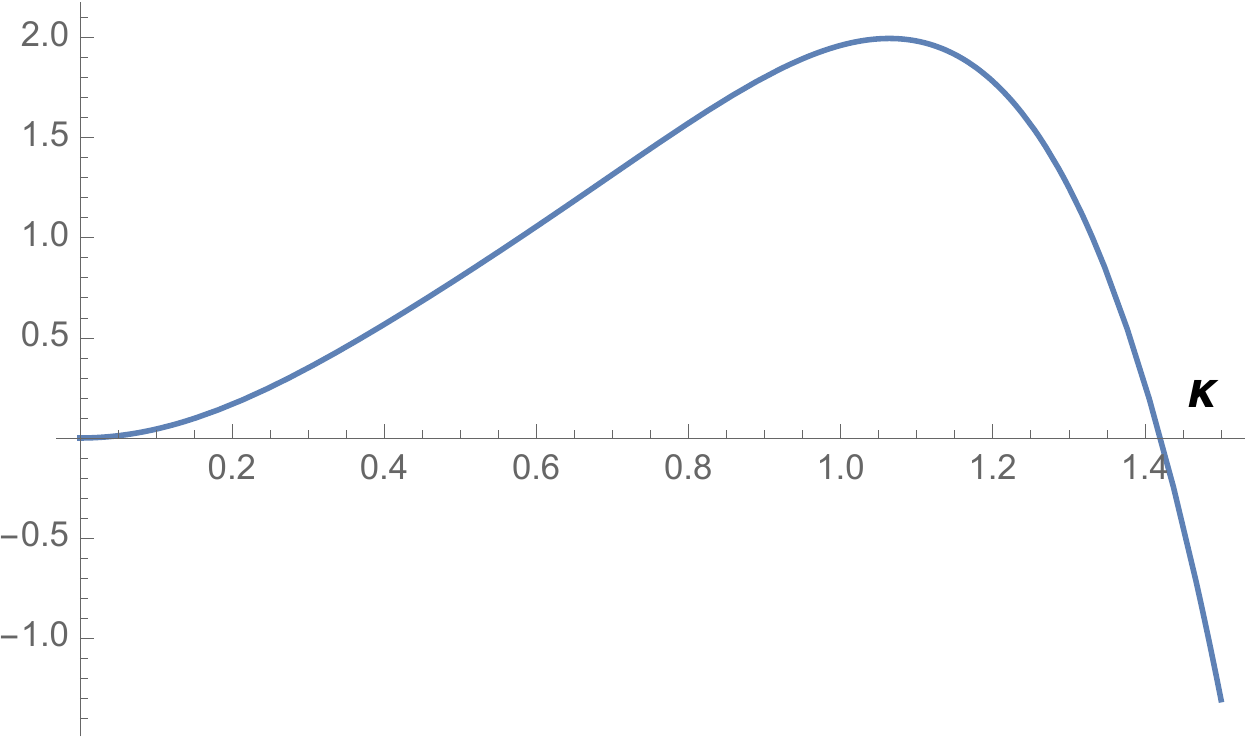}\quad
\caption{The graph of $i_4(\k)$ for $\k\in(0,1.5)$.}\label{fig:i4}
\end{figure}

We use \eqref{def:i4} and make an explicit calculation to show that 
\[
\lim_{\k\to 0+}\frac{i_4(\k)}{\k^2}=\frac92\quad\text{and}\quad\lim_{\k\to\infty}i_4(\k)=- \infty.
\]
Hence $\Delta(\k;0)>0$ for $\k>0$ sufficiently small, implying the modulational stability, and it is negative for $\k>0$ sufficiently large, implying the modulational instability. The intermediate value theorem asserts a root of $i_4$. A numerical evaluation of \eqref{def:i4} reveals a unique root $\k_c$, say, of $i_4$ over the interval $(0,\infty)$ such that $i_4(\k)>0$ if $0<\k<\k_c$ and it is negative if $\k_c<\k<\infty$. Upon close inspection (see Figure~\ref{fig:i4}), moreover, $\k_c=1.420\dots$. 

Therefore, a sufficiently small and $2\pi/\k$ periodic traveling wave of \eqref{E:FDCH} and \eqref{def:c1} is modulationally unstable if $\k>\k_c$, where $\k_c=1.420\dots$ is a unique root of $i_4$ in \eqref{def:i4} over the interval $(0,\infty)$. It is modulationally stable if $0<\k<\k_c$. The result qualitatively agrees with the Benjamin-Feir instability of a Stokes wave. The critical wave number compares reasonably well with that in \cite{Whitham1967} and \cite{BM1995}. The critical wave number for the Whitham equation is $1.146\dots$ (see \cite{HJ2}). 

\subsection*{For $\mathbf{T>0}$}

For $T>1/3$, since $\c(\k;T)$ and $(\k\c(\k;T))'$ increase monotonically over the interval $(0,\infty)$ and since $(\k\c(\k;T))'$ does not possess an extremum (see Figure~\ref{fig:cT}b), $i_1$, $i_2$, $i_3$ do not vanish over the interval $(0,\infty)$. A numerical evaluation reveals that $i_4$ changes its sign once over the interval $(0,\infty)$. Together, a sufficiently small and periodic traveling wave of \eqref{E:FDCH} and \eqref{def:cT1} is modulationally unstable, provided that the wave number is greater than a critical value, and modulationally stable otherwise, similarly to the gravity wave setting. Moreover, a numerical evaluation reveals that the critical wave number $\k_c(T)$, say, satisfies 
\[
\lim_{T\to \infty}\k_c(T)\sqrt{T}\approx 1.283.
\]

For $0<T<1/3$, on the other hand, $(\k\c(\k;T))'$ achieves a unique minimum over the interval $(0,\infty)$. Moreover, $i_2$ and $i_3$ each takes one transverse root over the interval $(0,\infty)$ (see Figure~\ref{fig:cT}c). Hence, $i_1$ through $i_4$ each contributes to the change in the modulational stability and instability.

Figure~\ref{f:ST} illustrates in the $\k$ and $\k\sqrt{T}$ plane the regions of modulational stability and instability for a sufficiently small and periodic traveling wave of \eqref{E:FDCH} and \eqref{def:cT1}. Along Curve~1, $i_1=0$ and the group speed achieves an extremum at the wave number $\k$. Curve~2 is associated with $i_2=0$, along which the group speed coincides with the phase speed in the limit as $\k\to0$. In the deep water limit, as $\k\to\infty$ while $\k\sqrt{T}$ is fixed, it is asymptotic to $\k=\frac94\k^2T-\frac34$. Curve~3 is associated with $i_3=0$, along which the phase speeds of the fundamental mode and the second harmonic coincide. In the deep water limit, it is asymptotic to $k^2T=\frac12$. Moreover, along Curve~4, $i_4$ vanishes because of the resonance of the dispersion and nonlinear effects. The ``lower" branch of Curve~4 passes through $\k=1.420\dots$, the critical wave number for $T=0$. The ``upper" branch passes through $\k\sqrt{T}=1.283\dots$, the limit of strong surface tension.

\begin{figure}[h] 
\includegraphics[scale=0.5]{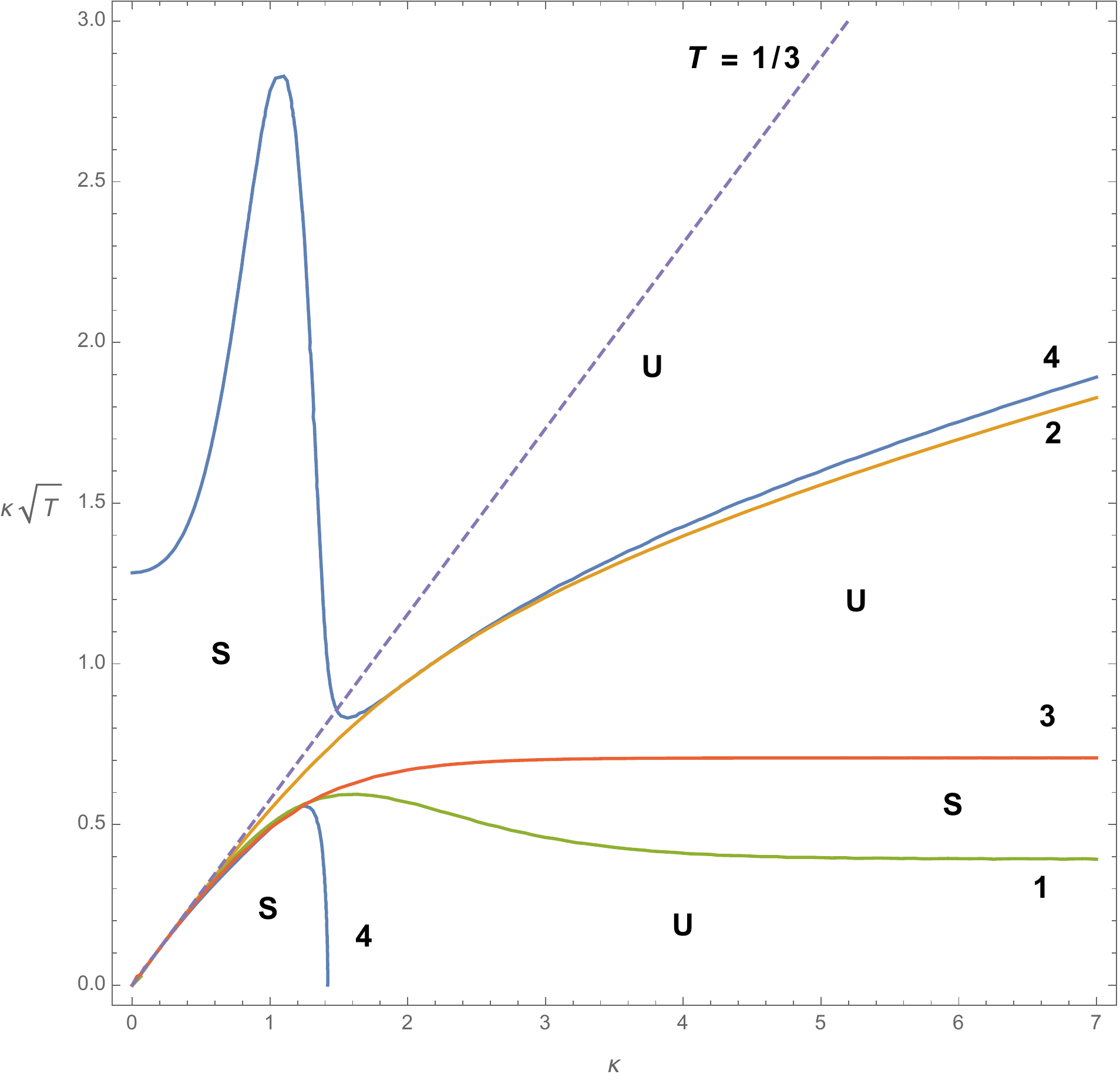}\quad
\caption{Stability diagram for sufficiently small and periodic traveling wave of \eqref{E:FDCH} and \eqref{def:cT1}. To interpret, for any $T>0$, one must envision a line through the origin with slope $T$. ``S" and ``U" denote stable and unstable regions. Solid curves labelled $1$ through $4$ represent the roots of $i_1$ through $i_4$ in \eqref{def:i1234}, respectively.}
\end{figure}\label{f:ST}

The result qualitatively agrees with those in \cite{Kawahara} and \cite{DR}, for instance, from formal asymptotic expansions for the physical problem, and it improves upon that in \cite{HJ3} for the Whitham equation, for which $\k_c(T)\sqrt{T}\to \infty$ as $T\to \infty$.

\section{Discussion}\label{sec:discussion}

\subsection*{The Camassa-Holm equation}

We may write \eqref{E:CH1}, in the case of $a=1/12$, after normalization of parameters, as 
\begin{equation}\label{E:CH}
\eta_t +c_\text{\Tiny CH}(|\partial_x|)\Big(2(1+\eta)^{3/2}-2\eta+\frac{7}{24}\eta \eta_{xx}+\frac{7}{48}\eta_x^2\Big)_x = 0,
\end{equation}
where
\begin{equation}\label{def:c-CH}
\widehat{c_\text{\Tiny CH}(|\partial_x|)f}(\k)=\frac{12-\k^2}{12+\k^2}\widehat{f}(\k).
\end{equation}
Note that $c_\text{\Tiny CH}(\k)$ approximates \eqref{def:c1}.

For any $\k>0$, we may repeat the argument in Section~\ref{sec:periodic} to determine sufficiently small and $2\pi/\k$ periodic traveling waves. Specifically, a two parameter family of periodic traveling waves of \eqref{E:CH} and \eqref{def:c-CH} exists, denoted 
\[
\eta(a,b;\k)(z)\qquad\text{where $z=\k(x-c(a,b;\k)t)$},
\]
for $a,b\in\mathbb{R}$ and $|a|,|b|$ sufficiently small; $\eta$ is $2\pi$ periodic, even, and smooth in $z$. Moreover,
\begin{align*}
\eta(a,b;\k)(z)=&b(1-c_\text{\Tiny CH}(\k))+a \cos z  + a^2 (h_0 + h_2 \cos 2z)+O(a(a^2+b^2)),\\
c(a,b;\k)=&c_\text{\Tiny CH}(\k)+b\Big(\frac32 - \frac{7}{24} \k^2 \Big)c_\text{\Tiny CH}(\k)(1-c_\text{\Tiny CH}(\k))+a^2 c_2+O(a(a^2+b^2)) 
\end{align*}
as $a,b \to 0$, where
\begin{gather*} 
h_0 =\frac{36-7 \k^2}{96(c_\text{\Tiny CH}(\k)-1)}, \qquad h_2 = \frac{(12-7 \k^2)c_\text{\Tiny CH}(2\k)}{32(c_\text{\Tiny CH}(\k)-c_\text{\Tiny CH}(2\k))},\\
c_2=c_\text{\Tiny CH}(\k)\Big( \frac{36-7 \k^2}{24} h_0+\frac{12-7 \k^2}{16}h_2-\frac{3}{32}\Big).
\end{gather*}
We then proceed as in Section~\ref{sec:MI}, to determine a modulational instability index
\[
\Delta_\text{\Tiny{CH}}(\k):=\frac{i_1(\k)i_2(\k)}{i_3(\k)}i_4(\k),
\]
where
\begin{align*}
i_1(\k)=&(\k c_\text{\Tiny{CH}}(\k))'',\\
i_2(\k)=&(\k c_\text{\Tiny{CH}}(\k))'-1, \\
i_3(\k)=&c_\text{\Tiny{CH}}(\k)-c_\text{\Tiny{CH}}(2\k),\\
i_4(\k)=&1296(c_\text{\Tiny{CH}}(2\k) i_2(\k)+2i_3(\k))-432 i_2(\k)i_3(\k) \\
&-1512 \k^2 c_\text{\Tiny{CH}}(2\k) i_2(\k)+49\k^4(9c_\text{\Tiny{CH}}(2\k) i_2(\k)-2i_3(\k)).
\end{align*}
We omit the details.
 
 A straightforward calculation shows that $\frac{i_2i_4}{i_3}(\k)< 0$ for any $\k > 0$ while $i_1(\k)$ changes its sign from negative to positive across $\k = 6$. Therefore, a sufficiently small and $2\pi/\k$ periodic traveling wave of \eqref{E:CH} and \eqref{def:c-CH} is modulationally unstable if $\k>6$. For other values of $a$ in \eqref{E:CH1}, the result is qualitatively the same. Thus, the Camassa-Holm equation seems to predict the Benjamin-Feir instability of a Stokes wave. But Resonance (R1) following Theorem~\ref{thm:FDCH} results in the instability in \eqref{E:CH} and \eqref{def:c-CH}, whereas it does not take place in the water wave problem.

\subsection*{The velocity equation}

The FDCH equation for the average horizontal velocity, after normalization of parameters, 
\begin{equation}\label{E:FDCHu}
u_t+c_{\rm ww}(|\partial_x|)u_x+\frac32 uu_x=-\Big(\frac{5}{12}uu_{xxx}+\frac{23}{24}u_xu_{xx}\Big),
\end{equation}
where $\c(|\partial_x|)$ is in \eqref{def:cT1}, differs from \eqref{E:FDCH} by higher power nonlinearities. We may repeat the arguments in Section~\ref{sec:periodic} and Section~\ref{sec:MI} to derive the modulational instability index, where \eqref{def:i4} is replaced by 
\begin{align} \label{def:iFDCHu}
i_u(\k;T)=&i_2(\k;T)+2i_3(\k;T) \\ 
&+\frac{1}{36}\k^2 (57i_2+34 i_3)(\k;T)+\frac{1}{324}\k^4(198 i_2+35i_3)(\k;T). \notag
\end{align}
We omit the details. The index formula for \eqref{E:FDCHu} and \eqref{def:cT1} agrees with that for the Whitham equation (see \cite{HJ2,HJ3} for details) except the terms in \eqref{def:iFDCHu} explicitly depending on $\k^2$ and $\k^4$, which higher derivative nonlinearities of \eqref{E:FDCHu} seem to contribute. 

For $T=0$, it is straightforward to show that a sufficiently small and $2\pi/\k$ periodic traveling wave of \eqref{E:FDCHu} and \eqref{def:cT1} is modulationally unstable if $\k>0.637\dots$. Thus the FDCH equation for the average horizontal velocity predicts the Benjamin-Feir instability of a Stokes wave. For $T>0$, the modulational instability result for \eqref{E:FDCHu} and \eqref{def:cT1} qualitatively agrees with that in \cite{HJ3} for the Whitham equation; see Figure~\ref{f:STu}. In particular, $\k_c(T)\sqrt{T}\to \infty$ as $T\to\infty$, where $\k_c(T)$ is a critical wave number, depending on $T$, whereas the limit is finite in the FDCH equation for the fluid surface displacement and the water wave problem. That means, the higher power nonlinearities of \eqref{E:FDCH} improve the modulational instability result in the presence of the effects of surface tension, not the higher derivative nonlinearities.

\begin{figure}[h] 
\includegraphics[scale=0.5]{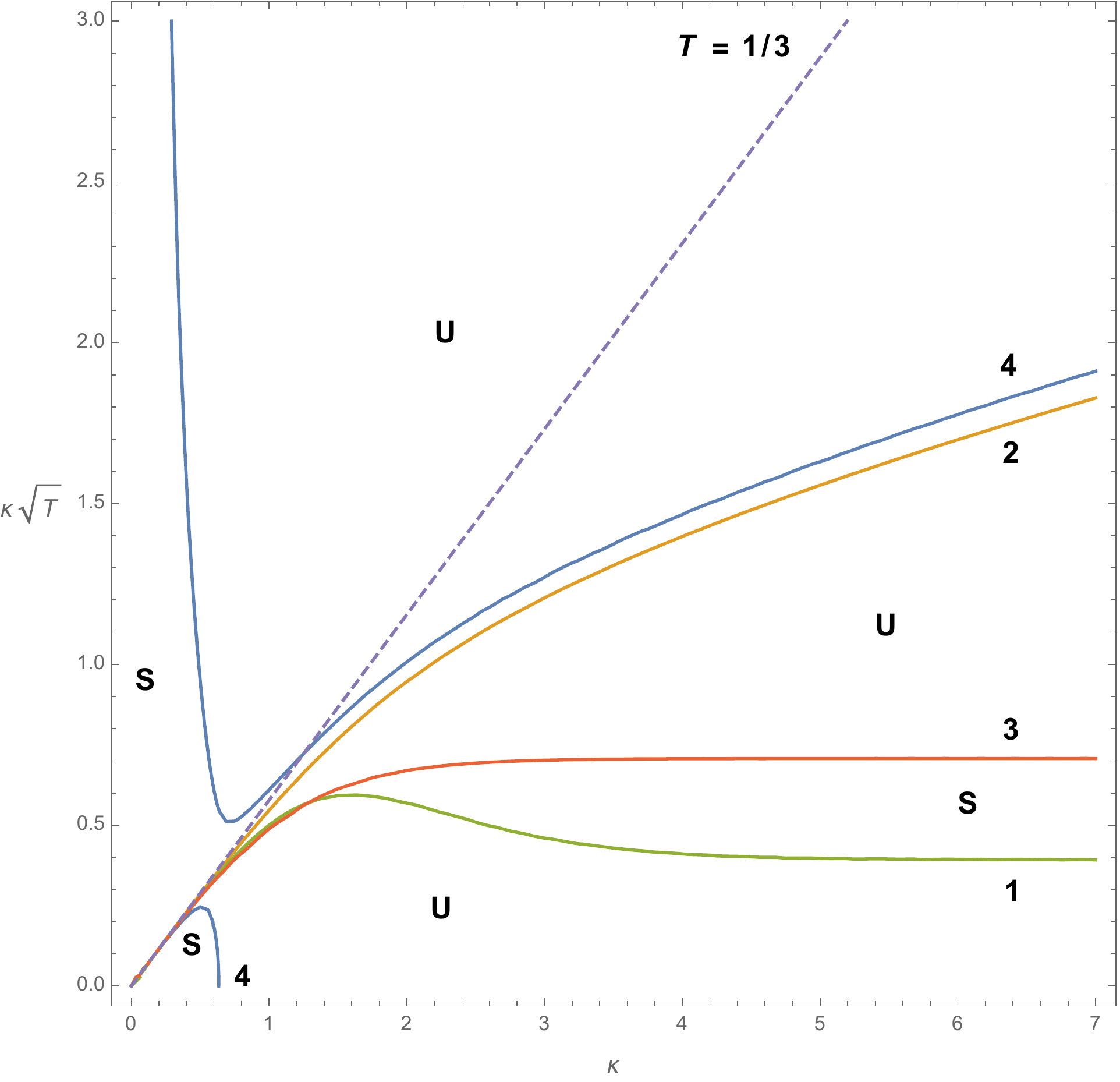}\quad
\caption{Stability diagram for sufficiently small and periodic wave trains of \eqref{E:FDCHu} and \eqref{def:cT1}.}
\end{figure}\label{f:STu}

\subsection*{Acknowledgements}
VMH is supported by the National Science Foundation's Faculty Early Career Development (CAREER) Award DMS-1352597, a Simons Fellowship in Mathematics, and the University of Illinois at Urbana-Champaign under the Arnold O. Beckman Research Award RB16227. She is grateful to the Department of Mathematics at Brown University for its generous hospitality.

\bibliographystyle{amsalpha}
\bibliography{stability.bib}
\end{document}